\documentclass[reqno]{tran-l}
\usepackage{graphicx}

 \newtheorem{thm}{Theorem}
 \newtheorem{cor}[thm]{Corollary}
 \newtheorem{lemma}[thm]{Lemma}
 \newtheorem{prop}[thm]{Proposition}
 \theoremstyle{definition}
 
 \newtheorem{ex}{Example}
 \theoremstyle{remark}

\newcommand{\C}{\mathbb{C}}
\newcommand{\D}{\mathbb{D}}

\newcommand{\calC}{\mathcal{C}}

\newcommand{\ovl}{\overline}

\newcommand{\rarrow}{\rightarrow}
\newcommand{\cphi}{C_{\phi}}
\newcommand{\phiinf}{\phi(\infty)}

\DeclareMathOperator{\RE}{Re}

\renewcommand{\phi}{\varphi}

\begin{document}

\title[Adjoints of composition operators]{Adjoints of composition operators \\ with rational symbol}

\author[Hammond, Moorhouse, and Robbins]{Christopher Hammond, Jennifer Moorhouse, and Marian E. Robbins}

\address{Department of Mathematics and Computer Science, Connecticut College,
Box 5384, 270 Mohegan Avenue, New London, CT 06320, U.S.A.}

\email{cnham@conncoll.edu}

\address{Department of Mathematics, Colgate University, 13 Oak
Drive, Hamilton, NY 13346, U.S.A.}

\email{jmoorhouse@colgate.edu}

\address{Mathematics Department, California Polytechnic State University, San Luis Obispo,
CA 93407, U.S.A.}

\email{mrobbins@calpoly.edu}


\date{May 3, 2007}

\dedicatory{}

\commby{}

%%% ----------------------------------------------------------------------
\maketitle

\thispagestyle{empty}

\begin{abstract}
Building on techniques developed by Cowen and
Gallardo-Guti\'{e}rrez, we find a concrete formula for the adjoint
of a composition operator with rational symbol acting on the Hardy
space $H^{2}$.  We consider some specific examples, comparing our
formula with several results that were previously known.
\end{abstract}

\section{Preliminaries}

Let $\D$ denote the open unit disk in the complex plane.  The
\textit{Hardy space} $H^{2}$ is the Hilbert space consisting of all
analytic functions $f(z)=\sum a_{n}z^{n}$ on $\D$ such that
\[
\left\|f\right\|_{2}=\sqrt{\sum_{n=0}^{\infty}|a_{n}|^{2}}<\infty\text{.}
\]
If $f(z)=\sum a_{n}z^{n}$ and $g(z)=\sum b_{n}z^{n}$ belong to
$H^{2}$, the inner product $\langle f,g\rangle$ can be written in
several ways.  For example,
\[
\langle f,g\rangle=\sum_{n=0}^{\infty}a_{n}\overline{b_{n}}
=\lim_{r\uparrow
1}\int_{0}^{2\pi}f(re^{i\theta})\overline{g(re^{i\theta})}\,\frac{d\theta}{2\pi}\text{.}
\]
Any function $f$ in $H^{2}$ can be extended to the boundary of $\D$
by means of radial limits; in particular, $f(\zeta)=\lim_{r\uparrow
1}f(r\zeta)$ exists for almost all $\zeta$ in $\partial\D$. (See
Theorem 2.2 in \cite{d}). Furthermore, we can write
\[
\langle
f,g\rangle=\int_{0}^{2\pi}f(e^{i\theta})\overline{g(e^{i\theta})}\,\frac{d\theta}
{2\pi}=\frac{1}{2\pi
i}\int_{\partial\D}f(\zeta)\overline{g(\zeta)}\,\frac{d\zeta}{\zeta}\text{.}
\]
It is often helpful to view $H^{2}$ as being a subspace of
$L^{2}(\partial\D)$.  Taking the basis
$\{z^{n}\}_{n=-\infty}^{\infty}$ for $L^{2}(\partial\D)$, we can
identify the Hardy space with the collection of functions whose
Fourier coefficients vanish for $n\leq-1$.

One important property of $H^{2}$ is that it is a
\textit{reproducing kernel Hilbert space}.  In other words, for any
point $w$ in $\D$ there is some function $K_{w}$ in $H^{2}$ (known
as a \textit{reproducing kernel function}) such that $\langle f,
K_{w}\rangle=f(w)$ for all $f$ in $H^{2}$. In the case of the Hardy
space, it is easy to see that $K_{w}(z)=1/(1-\overline{w}z)$.

At this point, we will introduce our principal object of study. Let
$\phi$ be an analytic map that takes $\D$ into itself. The
\textit{composition operator} $C_{\phi}$ on $H^{2}$ is defined by
the rule
\[
C_{\phi}(f)=f\circ\phi\text{.}
\]
It follows from Littlewood's Subordination Theorem (see Theorem 2.22
in \cite{cm}) that every such operator takes $H^{2}$ into itself.
These operators have received a good deal of attention in recent
years. Both \cite{cm} and \cite{s} provide an overview of many of
the results that are known.

\section{Adjoints}
One of the most fundamental questions relating to composition
operators is how to obtain a reasonable representation for their
adjoints. It is difficult to find a useful description for
$C_{\phi}^{\ast}$, apart from the elementary identity
\begin{equation}\label{kernel}
C_{\phi}^{\ast}(K_{w})=K_{\phi(w)}\text{.}
\end{equation}
(See Theorem 1.4 in \cite{cm}.)  In 1988, Cowen \cite{c} used this
fact to establish the first major result pertaining to the adjoints
of composition operators:

\begin{thm}[Cowen]\label{lfm}
Let
\[
\phi(z)=\frac{az+b}{cz+d}
\]
be a nonconstant linear fractional map that takes $\D$ into itself.
The adjoint $C_{\phi}^{\ast}$ can be written
$T_{g}C_{\sigma}T_{h}^{\ast}$, for
\[
g(z)=\frac{1}{-\overline{b}z+\overline{d}}\text{, }\sigma(z)
=\frac{\overline{a}z-\overline{c}}
{-\overline{b}z+\overline{d}}\text{, and }h(z)=cz+d\text{,}
\]
where $T_{g}$ and $T_{h}$ denote the Toeplitz operators with symbols
$g$ and $h$ respectively.
\end{thm}

\noindent While Cowen only stated this result for nonconstant
$\phi$, it is easy to see that it also holds when $\phi$ is
constant.  (In that case, $C_{\phi}$ and $C_{\sigma}$ can be
considered simply point-evaluation functionals.) It is sometimes
helpful to have a more concrete version of Cowen's adjoint formula.
Recalling that $T_{z}^{\ast}$ is the backward shift on $H^{2}$, we
see that
\begin{align}
(C_{\phi}^{\ast}f)(z)&=\left(\frac{1}{-\overline{b}z+\overline{d}}\right)
\left(\overline{c}
\left(\frac{f(\sigma(z))-f(0)}{\sigma(z)}\right)
+\overline{d}f(\sigma(z))\right)\nonumber\\
&=\left(\frac{1}{-\overline{b}z+\overline{d}}\right)
\left(\left(\frac{\overline{c}+\overline{d}\sigma(z)}{\sigma(z)}\right)f(\sigma(z))
-\frac{\overline{c}f(0)}{\sigma(z)}\right)\nonumber\\
&=\frac{(\overline{ad}-\overline{bc})z}{(\overline{a}z-\overline{c})(-\overline{b}z+\overline{d})}
f(\sigma(z))+\frac{\overline{c}f(0)}{\overline{c}-\overline{a}z}\text{.}\label{formula}
\end{align}
A similar calculation appears in \cite{h}.

In recent years, numerous authors have made the observation that
\begin{equation}\label{integral}
 (C_{\phi}^{\ast}f)(w)=\langle
f,K_{w}\circ\phi\rangle=\int_{0}^{2\pi}\frac{f(e^{i\theta})}
{1-\overline{\phi(e^{i\theta})}w}\frac{d\theta}{2\pi}\text{.}
\end{equation}
This fact seems particularly helpful when considering composition
operators induced by rational maps.  In an unpublished manuscript,
Bourdon \cite{b} uses it to find a representation for
$C_{\phi}^{\ast}$ when $\phi$ belongs to a certain class of
``quadratic fractional" maps. It is the principal tool used by
Effinger-Dean, Johnson, Reed, and Shapiro \cite{ejrs} to calculate
$\|C_{\phi}\|$ when $\phi$ is a rational map satisfying a particular
finiteness condition. Equation (\ref{integral}) is also the starting
point from which both Mart\'{i}n and Vukoti\'{c} \cite{mv} and Cowen
and Gallardo-Guti\'{e}rrez \cite{cg} attempt to describe the
adjoints of all composition operators with rational symbol. It is
the content of this last paper that provides the starting point for
our current discussion.

The results of Cowen and Gallardo-Guti\'{e}rrez are stated in terms
of \textit{multiple-valued weighted composition operators}. Suppose
that $\psi$ and $\sigma$ are a compatible pair of multiple-valued
analytic maps on $\D$ (in a sense the authors describe in their
paper), with $\sigma(\D)\subseteq\D$.  The operator
$W_{\psi,\sigma}$ is defined by the rule
\[
(W_{\psi,\sigma}f)(z)=\sum\psi(z)f(\sigma(z))\text{,}
\]
the sum being taken over all branches of the pair $\psi$ and
$\sigma$. Whenever we encounter such an operator in this paper, the
function $\psi$ will actually be defined in terms of~$\sigma$.

Before considering their adjoint theorem, we need to remind the
reader of a particular piece of notation.  If $f$ is a (possibly
multiple-valued) function on $\D$, we define the function
$\widetilde{f}$ on $\{z\in\C:|z|>1\}$ by the rule
\begin{equation}\label{tilde}
\widetilde{f}(z)=\overline{f\!\left(\frac{1}{\overline{z}}\right)}\text{.}
\end{equation}
Cowen and Gallardo-Guti\'{e}rrez state the following result:

\begin{thm}[Cowen and Gallardo-Guti\'{e}rrez]\label{main}
Let $\phi$ be a nonconstant rational map that takes $\D$ into
itself. The adjoint $C_{\phi}^{\ast}$ can be written
$BW_{\psi,\sigma}$, where $B$ denotes the backward shift operator
and $W_{\psi,\sigma}$ is the multiple-valued weighted composition
operator induced by $\sigma=1/\widetilde{\phi^{-1}}$ and
$\psi=\widetilde{(\phi^{-1})^{\prime}}/\widetilde{\phi^{-1}}$.
\end{thm}

\noindent Note that the function $\widetilde{(\phi^{-1})^{\prime}}$
in the numerator of $\psi$ represents the ``tilde transform" of
$(\phi^{-1})^{\prime}$, as defined in line (\ref{tilde}), rather
than the derivative of $\widetilde{\phi^{-1}}$.

As we shall see, Theorem \ref{main} is not entirely correct in all
cases.  We will begin by considering whether this result is valid
for linear fractional maps.

\section{Linear fractional examples}

If Theorem \ref{main} were to hold in general, it would certainly
have to agree with Theorem \ref{lfm} in the case of linear
fractional maps.  We shall show that these two theorems rarely yield
the same result. Let
\[
\phi(z)=\frac{az+b}{cz+d}
\]
be a nonconstant map that takes $\D$ into itself.  Note that
\[
\phi^{-1}(z)=\frac{dz-b}{-cz+a}\text{.}
\]
Using the notation of Theorem \ref{main}, we can write
\[
\widetilde{\phi^{-1}}(z)=\overline{\phi^{-1}(1/\overline{z})}
=\frac{-\overline{b}z+\overline{d}}{\overline{a}z-\overline{c}}\text{,}
\]
and thus
\[
\sigma(z)=1/\widetilde{\phi^{-1}}(z)
=\frac{\overline{a}z-\overline{c}}{-\overline{b}z+\overline{d}}\text{.}
\]
Also note that
\[
\widetilde{(\phi^{-1})^{\prime}}(z)
=\overline{\left(\frac{ad-bc}{\left(-\frac{c}{\overline{z}}+a\right)^{2}}\right)}
=\frac{(\overline{ad}-\overline{bc})z^{2}}{(\overline{a}z-\overline{c})^{2}}\text{,}
\]
from which it follows that
\[
\psi(z)=\widetilde{(\phi^{-1})^{\prime}}(z)/\widetilde{\phi^{-1}}(z)
=\frac{(\overline{ad}-\overline{bc})z^{2}}{(\overline{a}z-\overline{c})(-\overline{b}z+\overline{d})}\text{.}
\]
According to Theorem \ref{main}, the operator $C_{\phi}^{\ast}$ can
be written as the product $BW_{\psi,\sigma}$.  This would allow us
to write
\begin{equation}\label{noway}
(C_{\phi}^{\ast}f)(z)=\frac{\psi(z)f(\sigma(z))-\psi(0)f(\sigma(0))}{z}
=\frac{(\overline{ad}-\overline{bc})z}{(\overline{a}z-\overline{c})(-\overline{b}z+\overline{d})}
f(\sigma(z))\text{.}
\end{equation}
It is clear that line (\ref{noway}) agrees with line (\ref{formula})
only in one special case:  when $c=0$.  In other words, Theorem
\ref{main} is only valid for those linear fractional maps that are
actually linear.

There are several reasons that Theorem \ref{main} fails in this
situation. We will present a pair of examples, each demonstrating
one possible problem with the theorem.

\begin{ex}\label{ex1}
Consider the map
\[
\phi(z)=\frac{2z}{z+4}\text{,}
\]
which certainly takes $\D$ into itself.  Theorem \ref{main}, as
manifested in line (\ref{noway}), says that
\begin{equation}\label{noway2}
(C_{\phi}^{\ast}f)(z)
=\left(\frac{2z}{2z-1}\right)f(\sigma(z))\text{.}
\end{equation}
This statement cannot be correct, since this operator does not take
$H^{2}$ into itself. In particular, if $f(\sigma(1/2))=f(0)\neq 0$,
expression (\ref{noway2}) fails to be analytic in $\D$.

While we will provide a more detailed discussion at a later point,
the proper way to correct this difficulty is to follow the backward
shift with the orthogonal projection from $L^{2}(\partial\D)$ onto
$H^{2}$. In this particular example, we would obtain
\[
(C_{\phi}^{\ast}f)(z)=\left(\frac{2z}{2z-1}\right)f(\sigma(z))+\frac{f(0)}{1-2z}\text{.}
\]
This answer agrees with the adjoint formula given by line
(\ref{formula}).
\end{ex}

This projection is more than a notational nicety. The fact that the
projection is absent from the statement of Theorem \ref{main} causes
Corollary 3.9 in \cite{cg} to be incorrect.  In particular, the
kernel of $C_{\phi}^{\ast}$ includes not only those functions that
are sent to $0$ as a result of the backward shift, but also those
that are sent to $0$ due to the projection.

The next example illustrates a different problem with Theorem
\ref{main}.

\begin{ex}\label{ex2}
Let
\[
\phi(z)=\frac{z}{2z+4}\text{.}
\]
According to Theorem 2, as explained at the start of this section,
the adjoint should be given by the formula
\begin{equation}\label{noway3}
(C_{\phi}^{\ast}f)(z)=\left(\frac{z}{z-2}\right)f(\sigma(z))\text{.}
\end{equation}
Unlike the previous example, this expression always yields an
analytic function on $\D$.  Nevertheless, even if it did not
contradict Theorem \ref{lfm}, equation (\ref{noway3}) could not
possibly provide a correct formula for $C_{\phi}^{\ast}$. Since
$\varphi(0)=0$, line (\ref{kernel}) tells us that $C_{\phi}^{\ast}$
must fix the kernel function $K_{0}(z)=1$. It is obvious that the
operator defined above does not fix constant functions.

The difficulty here is more subtle -- and more fundamental -- than
that of the previous example.  The proof of Theorem \ref{main}
(i.e.\ Theorem 3.8 in \cite{cg}) relies on the assumption that
$\partial\phi^{-1}(\D)$, or at least one component thereof, is a
simple closed curve enclosing $\partial\D$. In this instance,
$\partial\phi^{-1}(\D)$ is a circle of radius $4/3$ centered at
$z=-8/3$.  In particular, $\partial\phi^{-1}(\D)$ does not enclose
$\partial{\D}$.  %(See Figure 1.)
Consequently one cannot use Cauchy's theorem to show that
\[
\int_{\partial\D}f(\phi(\zeta))\widetilde{g}(\zeta)\,\frac{d\zeta}{\zeta}=
\int_{\partial\phi^{-1}(\D)}f(\phi(\zeta))\widetilde{g}(\zeta)\,\frac{d\zeta}{\zeta}
\]
for polynomials $f$ and $g$, as asserted in the proof of Theorem
\ref{main}, since the integrand has a distinct singularity within
each circle (at $z=0$ and $z=-2$, respectively). Without this
purported equality, the remainder of the proof of Theorem \ref{main}
is invalid.

%\begin{figure}[ht]
%\scalebox{.65}{\includegraphics{shadegraph2a}}\caption{The set
%$\phi^{-1}(\D)$ for $\phi(z)=z/(2z+4)$, along with $\partial\D$.}
%\end{figure}
%Save as EPS; open with Photoshop (resolution 150); use color C0C0C0;
%save as EPS; (save as ASCII, no preview)

\end{ex}

Note that $\partial\phi^{-1}(\D)$ need not even be a circle.  For
example, if we take
\[
\phi(z)=\frac{z}{z+4}\text{,}
\]
then $\partial\phi^{-1}(\D)$ is the vertical line $\RE z=-2$.

While these examples are particular to linear fractional maps, the
problems they illustrate can manifest themselves with rational maps
of any degree.

\section{The main theorem}

From this point onward, we will assume that $\phi$ is a rational map
that takes $\D$ into itself.  Our goal is to apply certain aspects
of the proof of Theorem \ref{main} (Theorem 3.8 in \cite{cg}) to
obtain a general formula for $\cphi^{\ast}$. Before stating our main
results, we need to extract an important calculation from the proof
of Theorem 2.

\begin{lemma}[Cowen and Gallardo-Guti\'{e}rrez]\label{mainlem}
Let $\phi$ be a nonconstant rational map that takes $\D$ into
itself. Define the multiple-valued functions
\[
\sigma(z)=1/\widetilde{\phi^{-1}(z)}=1/\overline{\phi^{-1}(1/\overline{z})}
\]
and
\[
\psi(z)=\frac{z\sigma^{\prime}(z)}{\sigma(z)}\text{.}
\]
Suppose that $\partial\phi^{-1}(\D)$ is a bounded set.  If $f$ is a
rational function with no poles on $\partial\D$ and $g$ is a
polynomial, then
\begin{equation}\label{intrep}
\int_{\partial \phi^{-1}(\D)}f(\phi(\zeta))\widetilde{g}(\zeta)
\frac{d\zeta}{\zeta}=\int_{\partial\D}f(\zeta)\ovl{\sum\psi(\zeta)
g(\sigma(\zeta))}\frac{d\zeta}{\zeta}\text{.}
\end{equation}
In this expression, the function $\widetilde{g}$ is defined as in
line (\ref{tilde}) and the summation is taken over all branches of
$\sigma$.
\end{lemma}

While the map $\sigma$ here is the same as in Theorem \ref{main},
the function $\psi$ is slightly different. With a bit of work, the
$\psi$ from Theorem \ref{main} can be rewritten
$z^{2}\sigma^{\prime}(z)/\sigma(z)$. The reason for this discrepancy
is that here the function $\psi$ incorporates the action of the
backward shift.

It is important to recognize that the integral on the right-hand
side of (\ref{intrep}) does not necessarily represent an inner
product in $H^{2}$, or even in $L^{2}(\partial\D)$. As we have
already observed, the function $\sum\psi(g\circ\sigma)$ may have a
singularity in $\D$ and hence might not be analytic. Under certain
circumstances, it does not even belong to $L^{2}(\partial\D)$.

Before proceeding with the proof of this lemma, we need to consider
a minor technical detail. For $\phi$ analytic in $\D$, it is not
always the case that $\phi^{-1}(\partial\D)=\partial\phi^{-1}(\D)$.
These sets do coincide, however, whenever $\phi$ is analytic in a
neighborhood of $\overline{\phi^{-1}(\D)}$. Any rational map
$\phi:\D\rightarrow\D$ satisfies this condition, since it can only
have a finite number of poles, none of which can belong to
$\overline{\phi^{-1}(\D)}$.

\begin{proof}[Proof of Lemma \ref{mainlem}]
Fix a branch cut for $\phi^{-1}$ and consider all the branches
$\phi_{1}^{-1},$ $\phi_{2}^{-1},\ldots,\phi_{N}^{-1}$.  Define
$\Gamma_{j}=\phi_{j}^{-1}(\partial\D)$.  Note that the sets
$\Gamma_{j}$ are pairwise disjoint and that
\[
\bigcup_{j=1}^{N}\Gamma_{j}=\phi^{-1}(\partial\D)
=\partial\phi^{-1}(\D) \text{.}
\]
Thus
\begin{equation}\label{gams}
\int_{\partial \phi^{-1}(\D)}f(\phi(\zeta))\widetilde{g}(\zeta)
\frac{d\zeta}{\zeta}=\sum_{j=1}^N \int_{\Gamma_j} f(\phi(\zeta))
\widetilde{g}(\zeta) \frac{d\zeta}{\zeta}
\end{equation}
and we can apply the change of variables $\phi(\zeta)=\eta$ to each
of the integrals on the right-hand side of (\ref{gams}), since
$\phi$ is one-to-one on each $\Gamma_{j}$. Since
$\zeta=\phi_{j}^{-1}(\eta)$ for $\zeta$ in $\Gamma_{j}$, we obtain
\begin{align*}
\int_{\partial \phi^{-1}(\D)}f(\phi(\zeta))\widetilde{g}(\zeta)
\frac{d\zeta}{\zeta} & = \int_{\partial \D} f(\eta) \sum_{j=1}^N
\frac{\widetilde{g}(\phi_j^{-1}(\eta))}{\phi_j^{-1}(\eta)
\phi^{\prime}(\phi_j^{-1}(\eta))} d\eta\\
& =\int_{\partial \D} f(\eta) \sum_{j=1}^N
\frac{(\phi_j^{-1})^{\prime}(\eta)}{\phi_j^{-1}(\eta)}
\widetilde{g}(\phi_j^{-1}(\eta)) d\eta\text{.}
\end{align*}
If $\eta$ belongs to $\partial\D$, then
$\widetilde{g}(\phi_j^{-1}(\eta))=\ovl{g(\sigma_j(\eta))}$ and
$\phi_j^{-1}(\eta)=\ovl{\widetilde{\phi_j^{-1}}(\eta)}$.  One can
easily check that $(\widetilde{h})^{\prime}(z)=(-1/z^{2})
\widetilde{(h^{\prime})}(z)$ for all polynomials, and hence for all
functions that are analytic on a subset of $\D$. Thus a bit of
computation shows that
\[
\frac{\widetilde{(\phi_j^{-1})^{\prime}}(z)}{\widetilde{\phi_j^{-1}}(z)}
= \frac{z^2 \sigma_j^{\prime}(z)}{\sigma_j(z)}\text{.}
\]
Combining all of these observations, we have
\begin{align*}
\int_{\partial\phi^{-1}(\D)}f(\phi(\zeta))\widetilde{g}(\zeta)
\frac{d\zeta}{\zeta} &=\int_{\partial \D} f(\eta) \ovl{\sum_{j=1}^N
\frac{\eta^{2}\sigma_{j}^{\prime}(\eta)}{\sigma_{j}(\eta)}
g(\sigma_{j}(\eta))} d\eta\\
& = \int_{\partial \D} f(\eta) \ovl{\sum_{j=1}^N
\frac{\eta\sigma_{j}^{\prime}(\eta)}{\sigma_j(\eta)}
g(\sigma_{j}(\eta))} \frac{d\eta}{\eta}\text{,}
\end{align*}
as we had hoped to show.
\end{proof}

One consequence of this lemma is that, as long as it belongs to
$L^{2}(\partial\D)$, the function $\sum\psi(g\circ\sigma)$ is not
affected by the choice of branch cut.

While the multiple-valued function $\phi^{-1}$ can be quite
complicated in terms of its branch cuts, it can have only a single
pole. In particular, a rational map has only one limit as $|z|$ goes
to $\infty$. For the duration of this paper, we will write $\phiinf$
to denote this quantity; that is,
\[
\phiinf=\lim_{|z|\rightarrow\infty}\phi(z)\text{,}
\]
with the understanding that $\phiinf=\infty$ if no finite limit
exists.  The value of $\phiinf$ determines the basic geometric
properties of the set $\phi^{-1}(\D)$.  It is always the case that
$\phi^{-1}(\D)$ is an open set containing $\D$. As long as
$|\phiinf|\neq 1$, the boundary $\partial\phi^{-1}(\D)$ consists of
a finite number of bounded curves. If $|\phiinf|>1$, then
$\phi^{-1}(\D)$ is the union of the bounded regions enclosed within
$\partial\phi^{-1}(\D)$. If $|\phiinf|<1$, then $\phi^{-1}(\D)$ is
the unbounded exterior region.  If $|\phiinf|=1$, then at least one
component of $\partial\phi^{-1}(\D)$ must be unbounded, which means
that both $\phi^{-1}(\D)$ and its complement are unbounded sets.
%(See Figure 2.)

%\begin{figure}[ht]
%\scalebox{.65}{\includegraphics{shadegraph3a}}\caption{The set
%$\phi^{-1}(\D)$ for $\phi(z)=(z^{2}-1)/(z^{2}+z-4)$.}
%\end{figure}

For the purposes of this paper, we will need to consider each of
these cases separately.  We start with the most straightforward
situation.

\begin{prop}\label{prop1}
Suppose that $\phi$ is a nonconstant rational map that takes $\D$
into itself, with $|\phiinf|>1$.  Then
\[
(C_{\phi}^{\ast}f)(z)=\sum\psi(z)f(\sigma(z))+\frac{f(0)}
{1-\ovl{\phiinf}z}\text{,}
\]
where $\sigma$ and $\psi$ are the multiple-valued functions defined
in Lemma \ref{mainlem} and the summation is taken over all branches
of $\sigma$.
\end{prop}

\begin{proof}
Since the polynomials are dense in $H^{2}$, it suffices to consider
$\langle f,C_{\phi}^{\ast}(g)\rangle$ for polynomials $f$ and $g$.
Observe that $\widetilde{g}$, as defined in line (\ref{tilde}), is
analytic in $\{z:|z|>0\}$ and that $\widetilde{g}$ agrees with
$\overline{g}$ on the boundary of $\D$. Thus
\begin{equation}\label{booger}
\langle f,C_{\phi}^{\ast}(g)\rangle=\frac{1}{2\pi i}\int_{\partial
\D}
f(\phi(\zeta))\overline{g(\zeta)}\frac{d\zeta}{\zeta} \\
=\frac{1}{2\pi i}\int_{\partial\D} f(\phi(\zeta))
\widetilde{g}(\zeta) \frac{d\zeta}{\zeta}\text{.}
\end{equation}
We wish to transfer this last integral from $\partial\D$ to
$\partial\phi^{-1}(\D)$.

Since $|\phiinf|>1$, we know that $\phi^{-1}(\D)$ is a bounded set,
one component of which contains $\D$.  First of all, suppose $A$ is
a component of $\phi^{-1}(\D)$ that does not contain $\D$. Since
$\phi$ is analytic inside $A$ and the point $0$ does not belong to
$A$, we see that
\[
\int_{\partial
A}f(\phi(\zeta))\widetilde{g}(\zeta)\frac{d\zeta}{\zeta}=0\text{.}
\]
Now let $B$ denote the component of $\phi^{-1}(\D)$ that contains
$\D$. Since $\phi$ is analytic inside $B$, it is possible to deform
$\partial\D$ to $\partial B$ without passing over any singularities
of $\phi$ and without passing over 0. Thus we can apply Cauchy's
theorem and equation (\ref{booger}) to see that
\[
\langle f,\cphi^{\ast}(g)\rangle=\frac{1}{2\pi i}\int_{\partial B}
f(\phi(\zeta))\widetilde{g}(\zeta)\frac{d\zeta}{\zeta}
=\frac{1}{2\pi i}\int_{\partial\phi^{-1}(\D)}
f(\phi(\zeta))\widetilde{g}(\zeta)\frac{d\zeta}{\zeta}\text{.}
\]
Hence Lemma \ref{mainlem} dictates that
\begin{equation}\label{inner}
\langle f,\cphi^{\ast}(g)\rangle=\frac{1}{2\pi i}\int_{\partial\D}
f(\zeta)\ovl{\sum\psi(\zeta)g(\sigma(\zeta))}\frac{d\zeta}{\zeta}
\end{equation}
for any polynomials $f$ and $g$.

Consider the function $\sum\psi(g\circ\sigma)$ for a fixed
polynomial $g$.  We would like to show that $\sum\psi(g\circ\sigma)$
belongs to $L^{2}(\partial\D)$.  To that end, let us examine its
Fourier series with respect to the basis
$\{z^{n}\}_{n=-\infty}^{\infty}$. Since line (\ref{inner}) holds for
any polynomial $f$, we only need to compute the Fourier coefficients
corresponding to negative powers of $z$. Let $n$ be a natural number
and take $f(z)=z^{-n}$. Lemma \ref{mainlem} dictates that
\[
\frac{1}{2\pi i}\int_{\partial\D}
f(\zeta)\ovl{\sum\psi(\zeta)g(\sigma(\zeta))}\frac{d\zeta}{\zeta}=
\frac{1}{2\pi
i}\int_{\partial\phi^{-1}(\D)}\frac{\widetilde{g}(\zeta)}{\phi(\zeta)^{n}}
\frac{d\zeta}{\zeta}\text{.}
\]
Since $\phi$ is a rational map, the only singularities of the
integrand on the right-hand side occur when $\zeta=0$ or when
$\phi(\zeta)=0$; that is, within $\phi^{-1}(\D)$. Recall that the
set $\phi^{-1}(\D)$ is bounded. Therefore, for $R$ sufficiently
large, Cauchy's theorem shows that
\[
\frac{1}{2\pi
i}\int_{\partial\phi^{-1}(\D)}\frac{\widetilde{g}(\zeta)}{\phi(\zeta)^{n}}
\frac{d\zeta}{\zeta}= \frac{1}{2\pi
i}\int_{\mathcal{C}_{R}}\frac{\widetilde{g}(\zeta)}{\phi(\zeta)^{n}}
\frac{d\zeta}{\zeta}\text{,}
\]
where $\mathcal{C}_{R}$ denotes the circle of radius $R$ centered at
the origin.  Letting $R$ go to $\infty$, we see that
\[
\frac{1}{2\pi
i}\int_{\mathcal{C}_{R}}\frac{\widetilde{g}(\zeta)}{\phi(\zeta)^{n}}
\frac{d\zeta}{\zeta}=\frac{1}{2\pi
i}\int_{\partial\D}\frac{\widetilde{g}(R\zeta)}{\phi(R\zeta)^{n}}
\frac{d\zeta}{\zeta}\rightarrow\frac{\overline{g(0)}}{\phiinf^{n}}\text{.}
\]
In other words, the $(-n)$th Fourier coefficient of
$\sum\psi(g\circ\sigma)$ equals $g(0)/\overline{\phiinf^{n}}$, where
we understand this quantity to be $0$ if $\phiinf=\infty$.
Therefore, if $f(z)=\sum_{n=-N}^{N}a_{n}z^{n}$, we see that
\begin{align}\label{ineq}
\left|\frac{1}{2\pi i}\int_{\partial\D}
f(\zeta)\ovl{\sum\psi(\zeta)g(\sigma(\zeta))}\frac{d\zeta}{\zeta}\right|
&=\left|\sum_{n=1}^{N}a_{-n}\frac{\ovl{g(0)}}{\phiinf^{n}}
+\left\langle\sum_{n=0}^{N}a_{n}z^{n},C_{\phi}^{\ast}(g)
\right\rangle\right|\nonumber\\
&\leq\left\|f\right\|_{2}\left(\sqrt{\frac{|g(0)|^{2}}
{|\phiinf|^{2}-1}}+\left\|C_{\phi}\right\|\left\|g\right\|_{2}\right)\text{.}
\end{align}
Consequently $\sum\psi(g\circ\sigma)$ must belong to
$L^{2}(\partial\D)$.

In view of line (\ref{inner}), we see that $C_{\phi}^{\ast}(g)$ is
simply the projection of $\sum\psi(g\circ\sigma)$ onto $H^{2}$. Note
that
\[
\sum_{n=1}^{\infty}\frac{g(0)}{\overline{\phiinf^{n}}}\frac{1}{z^{n}}
=\frac{g(0)(\overline{\phiinf}z)^{-1}}{1-(\overline{\phiinf}z)^{-1}}
=\frac{g(0)}{\overline{\phiinf}z-1}\text{.}
\]
(For the purposes of this calculation, we can assume that $z$
belongs to $\partial\D$.)  Consequently the projection of
$\sum\psi(g\circ\sigma)$ onto $H^{2}$ equals
\[
\sum\psi(z)g(\sigma(z))+\frac{g(0)}{1-\overline{\phiinf}z}\text{,}
\]
which is precisely $C_{\phi}^{\ast}(g)$.
\end{proof}

Next we shall consider the situation when $\phi^{-1}$ has a pole
inside the disk; that is, when $|\phiinf|<1$.  The ideas underlying
our calculations will be similar to what we have already seen,
although $\phi^{-1}(\D)$ is an unbounded set and $\partial\D$ is not
enclosed within $\partial \phi^{-1}(\D)$. Hence we cannot use
Cauchy's theorem as in the previous proof, since in deforming
$\partial\D$ to $\partial \phi^{-1}(\D)$ we would pass over the
point $0$ and possibly some poles of $\phi$. Nevertheless, the
result we obtain looks quite familiar.

\begin{prop}\label{prop2}
Suppose that $\phi$ is a nonconstant rational map that takes $\D$
into itself, with $|\phiinf|<1$.  Then
\[
(C_{\phi}^{\ast}f)(z)=\sum\psi(z)f(\sigma(z))+\frac{f(0)}
{1-\ovl{\phiinf}z}\text{,}
\]
where $\sigma$ and $\psi$ are the multiple-valued functions defined
in Lemma \ref{mainlem} and the summation is taken over all branches
of $\sigma$.
\end{prop}

\begin{proof}
Let $f$ and $g$ be polynomials, and consider
\[
\langle f,C_{\phi}^{\ast}(g)\rangle=\frac{1}{2\pi
i}\int_{\partial\D}f(\phi(\zeta))
\widetilde{g}(\zeta)\frac{d\zeta}{\zeta}\text{.}
\]
Since $\partial\phi^{-1}(\D)$ is bounded, we can find $R$
sufficiently large so that $\mathcal{C}_{R}$, the circle of radius
$R$ centered at the origin, surrounds $\partial \phi^{-1}(\D)$. %(See Figure 3.)
Since $\partial\phi^{-1}(\D)$ encloses
$\phi^{-1}(\C\setminus\overline{\D})$, the map $\phi$ is analytic in
the region outside $\partial\phi^{-1}(\D)$.  Hence Cauchy's theorem
dictates that
\[
\int_{\calC_R} f(\phi(\zeta)) \widetilde{g}(\zeta)
 \frac{d\zeta}{\zeta} = \int_{\partial \D} f(\phi(\zeta)) \widetilde{g}(\zeta)
 \frac{d\zeta}{\zeta} + \int_{\partial \phi^{-1}(\D)} f(\phi(\zeta)) \widetilde{g}(\zeta)
 \frac{d\zeta}{\zeta}\text{,}
 \]
so that
\[
\langle f, \cphi^{\ast}(g)\rangle =-\frac{1}{2\pi
i}\int_{\partial\phi^{-1}(\D)} f(\phi(\zeta))\widetilde{g}(\zeta)
 \frac{d\zeta}{\zeta}+\frac{1}{2\pi i}\int_{\calC_R}f(\phi(\zeta))
\widetilde{g}(\zeta)
 \frac{d\zeta}{\zeta}\text{.}
 \]
Now we can apply the result of the Lemma \ref{mainlem},
understanding that the negative sign represents a reversal in the
orientation of $\partial\phi^{-1}(\D)$ and that the change of
variables results in a sign change. Therefore
\[
\langle f,\cphi^{\ast}(g)\rangle =\frac{1}{2\pi i}\int_{\partial\D}
f(\zeta)\overline{\sum\psi(\zeta)g(\sigma(\zeta))}\frac{d\zeta}{\zeta}
+\frac{1}{2\pi i}\int_{\calC_R}
f(\phi(\zeta))\widetilde{g}(\zeta)\frac{d\zeta}{\zeta}\text{.}
\]
If we let $R$ go to $\infty$, we see that
\[
\frac{1}{2\pi i}\int_{\calC_R} f(\phi(\zeta)) \widetilde{g}(\zeta)
\frac{d\zeta}{\zeta}=\frac{1}{2\pi i}\int_{\partial \D}
f(\phi(R\zeta)) \widetilde{g}(R\zeta)\frac{d\zeta}{\zeta} \rarrow
f(\phiinf) \ovl{g(0)}\text{.}
\]
(Recall our assumption that $|\phiinf|<1$.)  Notice that we can
rewrite
\[
f(\phiinf)\ovl{g(0)}=\langle C_{\phiinf}(f),C_0(g)\rangle\text{,}
\]
where $C_{\phiinf}$ and $C_{0}$ can simply be understood as
point-evaluation functionals. Theorem \ref{lfm} tells us that
$\langle C_{\phiinf}(f),C_0(g)\rangle=\langle f,
T_{\chi}C_{0}(g)\rangle$, where
\[
\chi(z)=\frac{1}{1-\ovl{\phiinf}z}\text{.}
\]
In other words,
\[
\langle f,\cphi^{\ast}(g)\rangle=\frac{1}{2\pi
i}\int_{\partial\D}f(\zeta)
\left(\overline{\sum\psi(\zeta)g(\sigma(\zeta))+\chi(\zeta)
g(0)}\right)\frac{d\zeta}{\zeta}
\]
for all polynomials $f$ and $g$.

%\begin{figure}[ht]
%\scalebox{.85}{\includegraphics{shadegraph4a}}\caption{The set
%$\phi^{-1}(\D)$ for $\phi(z)=(z^{2}+4)/(2z^{2}+4iz-8)$, along with
%$\partial\D$ and $\mathcal{C}_{R}$ for $R=5$.}
%\end{figure}

As in the proof of Proposition \ref{prop1}, let us consider the
Fourier series for the function $\sum\psi(g\circ\sigma)$. If
$f(z)=z^{-n}$, Lemma \ref{mainlem} dictates that
\[
\frac{1}{2\pi i}\int_{\partial\D}f(\zeta)
\overline{\sum\psi(\zeta)g(\sigma(\zeta))}\frac{d\zeta}{\zeta}=
\frac{1}{2\pi
i}\int_{\partial\phi^{-1}(\D)}\frac{\widetilde{g}(\zeta)}{\phi(\zeta)^{n}}
\frac{d\zeta}{\zeta}\text{.}
\]
The only singularities of the integrand on the right-hand side occur
within the set $\phi^{-1}(\D)$.  As we have already noted, the
curves that make up $\partial\phi^{-1}(\D)$ enclose the preimage of
the complement of the closed disk.  In other words, this integral
equals $0$ for all $n$, which means that the Fourier coefficients of
$\sum\psi(g\circ\sigma)$ corresponding to negative powers of $z$ are
all $0$.  Hence an argument similar to that of line (\ref{ineq})
shows that $\sum\psi(g\circ\sigma)$ belongs to $L^{2}(\partial\D)$,
and in fact to $H^{2}$. Therefore $C_{\phi}^{\ast}(g)$ is the
projection of $\sum\psi(g\circ\sigma)+\chi g(0)$ onto the Hardy
space. It is clear that $\chi$ belongs to $H^{2}$, so we obtain the
desired result.
\end{proof}

At first glance, it may appear that Propositions \ref{prop1} and
\ref{prop2} are saying the same thing, but in fact there is a subtle
difference.  As noted in the proof of Proposition \ref{prop2}, both
of the functions $\sum\psi(f\circ\sigma)$ and $\chi$ belong to
$H^{2}$. Thus, when $|\phiinf|<1$, we can write
\[
C_{\phi}^{\ast}=W_{\psi,\sigma}+T_{\chi}C_{0}\text{,}
\]
where all of these operators actually take $H^{2}$ into itself.

Finally we turn our attention to the case where $|\phiinf|=1$. In
this situation, as we have already mentioned, the geometric
properties of $\phi^{-1}(\D)$ can be rather complicated.
Furthermore, as we can deduce from the next proposition, the
function $\sum\psi(z)f(\sigma(z))$ does not belong to
$L^{2}(\partial\D)$. Rather than dealing with this case on its own
terms, we will obtain the desired result as a consequence of
Proposition \ref{prop2}.

\begin{prop}\label{prop3}
Suppose that $\phi$ is a nonconstant rational map that takes $\D$
into itself, with $|\phiinf|=1$.  Then
\[
(C_{\phi}^{\ast}f)(z)=\sum\psi(z)f(\sigma(z))+\frac{f(0)}
{1-\ovl{\phiinf}z}\text{,}
\]
where $\sigma$ and $\psi$ are the multiple-valued functions defined
in Lemma \ref{mainlem} and the summation is taken over all branches
of $\sigma$.
\end{prop}

\begin{proof}
Fix a real number $r$, with $0<r<1$, and consider the map
$r\phi(z)$. Observe that $|r\phiinf|<1$, so Proposition \ref{prop2}
applies to the operator $C_{r\phi}^{\ast}$.  Define $\sigma$ and
$\psi$ in terms of the original map $\phi$; we shall state the
formula for $C_{r\phi}^{\ast}$ with respect to these functions.
Since $(r\phi)^{-1}(z)=\phi^{-1}(z/r)$, it follows that
\[
\frac{1}{\overline{(r\phi)^{-1}(1/\overline{z})}}
=\frac{1}{\overline{\phi^{-1}(1/\overline{rz})}}=\sigma(rz)
\]
and
\[
\frac{z(\sigma(rz))^{\prime}}{\sigma(rz)}=
\frac{(rz)\sigma^{\prime}(rz)}{\sigma(rz)}=\psi(rz)\text{.}
\]
Consequently
\begin{equation}\label{radj}
(C_{r\phi}^{\ast}f)(z)=\sum\psi(rz)f(\sigma(rz))
+\frac{f(0)}{1-\overline{\phiinf}rz}\text{.}
\end{equation}
Observe that $C_{r\phi}=C_{\phi}C_{\rho}$, where $\rho(z)=rz$. Since
$C_{\rho}$ is self-adjoint (as we can see from Theorem \ref{lfm}),
it follows that
\[
(C_{r\phi}^{\ast}f)(z)=(C_{\rho}^{\ast}C_{\phi}^{\ast}f)(z)
=(C_{\rho}C_{\phi}^{\ast}f)(z)=(C_{\phi}^{\ast}f)(rz)\text{.}
\]
Combining this observation with line (\ref{radj}), we see that
\[
(C_{\phi}^{\ast}f)(rz)=\sum\psi(rz)f(\sigma(rz))
+\frac{f(0)}{1-\overline{\phiinf}rz}
\]
for all $z$ in $\D$.  Making the substitution $\zeta=rz$, we
conclude that
\[
(C_{\phi}^{\ast}f)(\zeta)=\sum\psi(\zeta)f(\sigma(\zeta))
+\frac{f(0)}{1-\overline{\phiinf}\zeta}
\]
whenever $|\zeta|<r$.  Since $r$ is arbitrary, our assertion
follows.
\end{proof}

Combining Propositions \ref{prop1}, \ref{prop2}, and \ref{prop3}, we
can now state our main result.

\begin{thm}\label{correct}
Suppose that $\phi$ is a rational map that takes $\D$ into itself.
Then
\[
(C_{\phi}^{\ast}f)(z)=\sum\psi(z)f(\sigma(z))+\frac{f(0)}
{1-\ovl{\phiinf}z}\text{,}
\]
where
\begin{align*}
\sigma(z)&=1/\widetilde{\phi^{-1}(z)}=1/\overline{\phi^{-1}(1/\overline{z})}\text{,}\\
\psi(z)&=\frac{z\sigma^{\prime}(z)}{\sigma(z)}\text{,}\\
\phiinf&=\lim_{|z|\rightarrow\infty}\phi(z)\text{,}
\end{align*}
and the summation is taken over all branches of $\sigma$.
\end{thm}

Notice that this result agrees exactly with line (\ref{formula}),
the adjoint formula for composition operators with linear fractional
symbol.  Furthermore, we obtain the following corollary.

\begin{cor}
Suppose that $\phi$ is a rational function that takes $\D$ into
itself.  If $\phiinf=\infty$, then $C_{\phi}^{\ast}$ is a
(multiple-valued) weighted composition operator.
\end{cor}

\section{Further examples}

Prior to this paper, there were a few concrete examples of rational
$\phi$ (in addition to the linear fractional maps) for which the
adjoint $C_{\phi}^{\ast}$ could be described precisely. We will
conclude by considering several such results in the context of
Theorem \ref{correct}.

\begin{ex}
Let $\phi(z)=z^{m}$ for some natural number $m$.  It is easy to
calculate $C_{\phi}^{\ast}$ simply by considering its action on the
orthonormal basis $\{z^{n}\}_{n=0}^{\infty}$.  (See Exercise 9.1.1
in \cite{cm}.) Applying Theorem \ref{correct}, we see that
\[
(C_{\phi}^{\ast}f)(z)=\sum_{j=1}^{m}\frac{f(\sigma_{j}(z))}{m}\text{,}
\]
where $\sigma_{1},\sigma_{2},\ldots,\sigma_{m}$ constitute all the
branches of the function $\sqrt[m]{z}$. This result agrees with the
formulas previously stated in \cite{mv} and \cite{m}.
\end{ex}

\begin{ex}
The one specific example discussed by Cowen and
Gallardo-Guti\'{e}rrez \cite{cg} is the map $\phi(z)=(z^{2}+z)/2$.
Let us consider a slightly more general case:
\[
\phi(z)=az^{2}+bz\text{,}
\]
with $a\neq 0$.  Note that $\phiinf=\infty$, so $C_{\phi}^{\ast}$ is
actually a multiple-valued weighted composition operator.  It is
easy to show that
\[
\sigma(z)=\frac{\overline{b}z\pm
\sqrt{\overline{b}^{2}z^{2}+4\overline{a}z}}{2}
\]
and
\[
\psi(z)=\frac{z\sigma^{\prime}(z)}{\sigma(z)}
=\frac{1}{2}\left(1\pm\frac{\overline{b}\sqrt{\overline{b}^{2}z^{2}+4\overline{a}z}}
{\overline{b}^{2}z+4\overline{a}}\right)\text{.}
\]
Therefore Theorem \ref{correct} says that
\[
(C_{\phi}^{\ast}f)(z)=\sum_{j=1}^{2}\frac{1}{2}\left(1
+(-1)^{j}\frac{\overline{b}\sqrt{\overline{b}^{2}z^{2}+4\overline{a}z}}
{\overline{b}^{2}z+4\overline{a}}\right)
f\!\left(\frac{\overline{b}z+
(-1)^{j}\sqrt{\overline{b}^{2}z^{2}+4\overline{a}z}}{2}\right)\text{.}
\]
When $a=b=1/2$, this expression agrees with the calculation found in
\cite{cg}.
\end{ex}

\begin{ex}
As we mentioned earlier, several years ago Bourdon \cite{b}
calculated the adjoint $C_{\phi}^{\ast}$ for $\phi$ belonging to a
particular class of rational maps. One example he considered was
\[
\phi(z)=\frac{z^{2}-6z+9}{z^{2}-10z+13}\text{.}
\]
Note that
\[
\sigma(z)=\frac{3z-5\pm 2\sqrt{3-2z}}{9z-13}
\]
and
\[
\psi(z)=\frac{\pm 2z}{\sqrt{3-2z}(3z-4\pm\sqrt{3-2z})}\text{.}
\]
Therefore $(C_{\phi}^{\ast}f)(z)$ equals
\[
\sum_{j=1}^{2}\frac{(-1)^{j}\,
2z}{\sqrt{3-2z}(3z-4+(-1)^{j}\sqrt{3-2z})}f\!\left(\frac{3z-5+(-1)^{j}\,
2\sqrt{3-2z}}{9z-13}\right)+\frac{f(0)}{1-z}\text{.}
\]
While the two formulas look somewhat different, it is not difficult
to show that this result is identical to Bourdon's.
\end{ex}

\section*{Acknowledgements}
The authors are indebted to Paul Bourdon for his many helpful
comments and suggestions.  The first and second authors also wish to
thank the Mathematics Department at California Polytechnic State
University for its generous hospitality during the winter quarter of
2007.


\begin{thebibliography}{99}

\bibitem{b} P. S. Bourdon, Some adjoint formulas for composition
operators on $H^{2}$, unpublished manuscript, 2002.

\bibitem{c}C. C. Cowen, Linear fractional composition operators on $H^{2}$,
\textit{Integral Equations Operator Theory}, \textbf{11} (1988),
151--160.

\bibitem{cg} C. C. Cowen and E. A. Gallardo-Guti\'{e}rrez,
A new class of operators and a description of adjoints of
composition operators, \textit{J. Funct. Anal.}, \textbf{238}
(2006), 447--462.

\bibitem{cm} C. C. Cowen and B. D. MacCluer,
\textit{Composition Operators on Spaces of Analytic Functions}, CRC
Press, Boca Raton, 1995.

\bibitem{d} P. L. Duren, \textit{Theory of $H^{p}$ spaces}, Academic
Press, New York, 1970.

\bibitem{ejrs} S. Effinger-Dean, A. Johnson, J. Reed, and J.
Shapiro, Norms of composition operators with rational symbol,
\textit{J. Math. Anal. Appl.}, \textbf{324} (2006) 1062--1072.

\bibitem{h} C. Hammond, On the norm of a composition operator
with linear fractional symbol, \textit{Acta Sci. Math. (Szeged)},
\textbf{69} (2003), 813--829.

\bibitem{mv} M. J. Mart\'{i}n and D. Vukoti\'{c}, Adjoints of
composition operators on Hilbert spaces of analytic functions,
\textit{J. Funct. Anal.}, \textbf{238} (2006), 298--312.

\bibitem{m} J. N. McDonald, Adjoints of a class of composition
operators, \textit{Proc. Amer. Math. Soc.}, \textbf{131} (2003),
601--606.

\bibitem{s} J. H. Shapiro, \textit{Composition Operators and
Classical Function Theory}, Springer-Verlag, New York, 1993.

\end{thebibliography}
\end{document}